\documentclass[12pt]{article}
\usepackage{amssymb}                         
\usepackage[T2A]{fontenc}
\usepackage{enumerate,amsmath,graphicx}
\usepackage{graphicx}

\begin{document}

\pagestyle{plain}
\title{Reducing the Number of Changeover Constraints in a MIP Formulation of
a Continuous-Time Scheduling Problem}

\author{Pavel A. Borisovsky$^1$, Anton V. Eremeev$^2$, Josef Kallrath$^3$}

\maketitle
\begin{center}
$\mbox{ }^1$ Omsk State University n.a. F.M. Dostoevskiy, \\
55a, pr. Mira, 644077, Omsk, Russia\\
$\mbox{ }^2$ Omsk Branch of Sobolev Institute of mathematics,\\
13 Pevtsov str., 644099, Omsk, Russia\\
$\mbox{ }^3$ BASF SE, 67056, Ludwigshafen, Germany
\end{center}

\begin{abstract}
In this paper, we develop a new formulation of changeover
constraints for mixed integer programming problem~(MIP) that
emerges in solving a short-term production scheduling problem. The
new model requires fewer constraints than the original formulation
and this often leads to shorter computation time of a MIP solver.
Besides that, the new formulation is more flexible if the time
windows for changeover tasks are given.\\

{\bf Keywords:} Short-term scheduling, mixed integer programming
model, changeover task

\end{abstract}

\section{Introduction}

Scheduling the processes of an industrial plant involves
a large number of objects such as processing units,
tasks, intermediate and final products (states)
and a set of complex relations between them.
The production scheduling problem basically consists in selection
of a set of tasks to be performed and construction of a
schedule complying with the technological requirements and
satisfying as much as possible the given demands
on a final production.

In earlier studies, a great deal of work was done concerning the
modelling of such scheduling problems in a form of Mixed Integer
Programs (MIP), see~\cite{Cho2002, JFKV2006, Floudas98,
PG1998,FloudasKTE2009}. Solving such a scheduling problem as mixed
integer program  meets serious difficulties when the number of
integer variables and constraints increases. To overcome this,
several variants of decomposition scheme were proposed, in which
the planning period is split into a sequence of smaller horizons,
and in each horizon a short-term MIP problem is solved separately.
The shorter the length of horizons, the less products are
considered in each of the short-term scheduling problems and the
fewer choices are available for optimization routine in each
horizon. Increasing the length of horizons usually leads to
improvement of solution quality but increases the CPU time and
memory requirements, which may become prohibitive at some point.

The  main objective of the paper is to reduce the number of
constraints in MIP formulation of one-horizon short-term
scheduling problem proposed in~\cite{FloudasKTE2009}. Reduced
number of constraints allows us to increase the size of manageable
subproblems and improve the overall performance of decomposition
algorithm~\cite{BEFSK}.

\section{Problem Statement} \label{sec:Form}

A modern chemical production plant is organized as a flexible
automated system that contains a number of multipurpose production
units and produces a number of products including the final
products and intermediate states. The control of the production
process involves two main problems: the first one is to choose the
most appropriate production plan (the set and the order of
reactions to be performed) that can satisfy the market
requirements on amount and assortment of the final product, and
the second one is to build an optimal production schedule for the
chosen plan. These two problems are closely related and for the
most effective control they can not be considered separately.

The input data for the basic production scheduling problem~
\cite{Floudas98,FloudasKTE2009} consists of:
\begin{itemize}
\item The set of production units.
\item The set of states, including raw materials, intermediate and
final product.
\item The set of tasks.
\item The suitability table assigning a suitable unit to each task.
Here we suppose without loss of generality that each task can be
performed only on one unit.
\item The minimal and maximal load of each unit when performing a task.
\item The production recipes represented in a form of a {\it State
Task Network (STN)}: for each task, a set of states being consumed
and produced are given, together with their consumption and
production coefficients.
\item The processing rates for each task, which is the amount of product produced by a task in one hour.
\item The demands and due dates for the final products.
\end{itemize}

The problem asks for a feasible schedule of tasks. The main
objectives are to minimize the underproduced amounts and
deviations from the due times for the final products.

In this paper, we consider an industrial plant of a special
structure. The principal production processes are performed on
the {\it main units}, the number of which is relatively small.
In addition there are subsidiary units that serve for transportation
and storage of the materials.
The main units need a changeover task when switching from one
production task to another. The durations of the changeover tasks
are sequence dependent, i.e. they are defined for the pairs of
production tasks and their order. The considered plant has special
requirements for the changeover tasks:

(i) No changeovers on main units are allowed between 4~a.m and
7~a.m or between 4~p.m. and 7~p.m.

(ii) There are two types of main tasks (the corresponding sets are
denoted by $I^{\rm p1}$ and $I^{\rm np1}$), that define two types
of changeovers. Let $I^{\rm c1}$ be the set of changeover tasks
from $i^1\in I^{\rm p1}$ to $i^2 \in I^{\rm np1}$. Such tasks may
be performed only between 7:30~a.m. and 1:30~p.m. The set of all
the rest of changeover tasks on the main units is denoted by
$I^{\rm c}$. So for any production unit there are at most two
changeover tasks available (one from $I^{\rm c1}$ and one from $I^{\rm c}$).


\section{Decomposition Approach} \label{sec:Decomp}
The decomposition approach was proposed in \cite{Floudas98} and in
\cite{FloudasKTE2009} it was adapted to the problem under
consideration. The whole planning period is split into smaller
horizons and a series of scheduling problems corresponding to
these horizons is solved sequentially. Two MIP models are used:
the {\it upper-level decomposition model (Level-1)} determines the
time partitioning and assigns the demands to small horizons; the
{\it lower-level short-term scheduling model (Level-2)} is used to
build a schedule in each horizon.

\subsection{Upper-Level Decomposition Model (Level-1)}

For the considered problem, the planning period (about one month)
is split into small horizons with the length of 12 hours. The aim
of Level-1 model is to determine how many small horizons will be
chosen for the short-term scheduling, and which tasks and states
will be considered. The complete model  will not be reproduced
here because it has many particular details that do not play an
important role for our study. Instead, we will only give a list of
the most essential conditions that are incorporated in the model:
\begin{itemize}
\item If some horizon is selected in the Level-1 model,
 then all preceding horizons are selected as well
(the alternation of selected and non selected horizons is not
allowed).
\item Product having due dates on some horizon must be included
into this horizon.
\item If some product is included then all intermediate
states corresponding to this product must be included.
\item The number of binary variables of the short-term scheduling
is bounded by a given tunable parameter.
\item The load times of the main units are bounded
by a given tunable parameter.
\end{itemize}

The objective includes the minimization of the number of selected
horizons, states,  maximization of the demand amount selected and
some secondary technical criteria (see the details in
\cite{FloudasKTE2009}). Despite the fact that the Level-1 problem
is of mixed integer programming type, the solving time of this
problem usually does not exceed several seconds.

\subsection{Lower-Level Short-Term Scheduling Model (Level-2)}

The MIP model for the short-term scheduling is based on the
continuous time concept~\cite{Floudas98}. The aim of the model is
to determine the set of operations to perform, their start and
finish times, and amounts of produced and consumed states for each
task.

The selection of tasks and their sequencing is done by introducing
the set of {\it event points} $N = \{1,2,...,N^{\max}\}$ for each
unit. An event point represents a relative position in the
sequence where a task can be assigned. Note that an event point
for some unit can be empty, i.e. have no tasks assigned. The
binary variable $wv(i,n)$ equals 1 iff a task $i$ is scheduled on
the corresponding unit at event point $n$ (recall that by
assumption each task has only one suitable unit). If a task is
scheduled at some event point, then we will say that the event
point has an active task. The material and timing conditions are
modelled by real-valued variables: $Ts(i,n)$ and $Tf(i,n)$ give
the start and the finish time of task~$i$ at event point~$n$, and
$B(i,n)$ defines the amount of processed material (in case
$wv(i,n)=0$ the values of these variables are unimportant).

The system of constraints includes the following conditions:

\begin{itemize}
\item For each unit an event point contains at most one task.
\item The amount of material processed by a task lies between given minimal and maximal bounds.
\item At each event point the produced and consumed amount
of some state must be balanced according to production
recipes.
\item The tasks producing and consuming some intermediate
state must be synchronized in time.
\item The tasks performing on the same unit must not overlap in time.
\item The execution time of a task is defined by its amount.
\item The main units may require changeover tasks,
i.e. if a task $j$ is assigned after task $i$ on the same main
unit, then a special changeover task $i^c$ must be scheduled
between them. For each pair of tasks $i$ and $j$, the changeover
duration is given.
\end{itemize}

The primary objective is minimization of the non-delivered amount
for final products with positive demands (underproduction). A set
of secondary objectives can be introduced: minimization of the
deviations from the delivery due dates, the overproduced amounts,
the number and duration of changeover and processing tasks, and so
on.

In this paper, we use the MIP model from \cite{FloudasKTE2009},
but reformulate the part concerning the changeover tasks. The new
model is equivalent to the original one, but has smaller number of
constraints. In what follows, we do not describe the whole set of
constraints of the Level-2 model, but only the part that was
modified and its new version.

\section{Former Model of Changeover Tasks} \label{sec:former_model}
The following notation
will be used below:
\begin{itemize}
\item $M$ is a set of main units.
\item $H$ is is the length of the planning horizon.
\item  $I_u$ is the set of tasks that can be performed on unit
$u$.
\item $Ctime_{i,i'}$ is a changeover time required to switch a unit from
task $i$ to task $i'$.
\item $I^{\rm p}$ is the set of the processing tasks (task of the main units
excluding  changeover tasks).
\end{itemize}

The binary variables $wv(i,n)$ equal 1 iff a task $i$ is active on
its suitable unit at event point $n$. Let us introduce binary
variable $x(i',i,n)$ that equals to 1 iff there is a changeover
from task~$i'$ occurring at event point~$n$ to some other task~$i$
($i \ne i'$) occurring at a later event point~$n'$ ($n' > n$) and
no other task is active between $n$ and $n'$ on the
unit suitable for $i'$ and $i$.\\

Let us reproduce the constraints (27a)-(28b) from
\cite{FloudasKTE2009}.\\

When task $i'$ is not active at an event point $n$, the variable
 $x(i',i,n)$ must be zero:
\begin{equation}
x(i',i,n) \le wv(i',n) \ \forall u\in M,\ i,i'\in I_u,i\neq i',
Ctime_{i',i}>0, \ n<N^{\max}.
\end{equation}

The changeover from task~$i'\in I_u$ to $i\in I_u$ is not
activated at event point~$n$ if there is a task~$i''\ne i$ at some
event point $n'>n$ and all event points in subsequence
$n+1,n+2,...,n'-1$ are empty on unit $u$:
$$
x(i',i,n) \le wv(i,n') + \left(1-\sum_{i''\in I^{\rm p} \cap
I_u}wv(i'',n') \right) +
$$
$$
\sum_{i''\in I^{\rm p} \cap I_u \mbox{\ }} \sum_{n'': n<n''<n'}
wv(i'',n'') \ \
$$
\begin{equation}
\label{eqn:ch2} \forall u\in M,\ i,i'\in I_u,i\neq i',
Ctime_{i',i}>0, \ n<N^{\max},\ n<n'\leq N^{\max}.
\end{equation}

The changeover from $i'$ to $i$ is activated if task $i'$ is
active in event point~$n$, task $i$ is active in event
point~$n'>n$, and there is no active task in event points
$n+1,n+2,...,n'-1$:

$$
x(i',i,n) \ge wv(i',n) + wv(i,n') - 1 -
 \sum_{i''\in I^{\rm p} \cap I_U\mbox{\ }} \sum_{n'': n<n''<n'} wv(i'',n'')
$$
\begin{equation}
\label{eqn:ch3} \forall u\in M,\ i,i'\in I_u,i\neq i',
Ctime_{i',i}>0, \ n<N^{\max},\ n<n' \leq N^{\max}.
\end{equation}

In the following two equations, the changeover task is allocated
in event point $n+1$ if there are non-zero values of $x(i',i,n)$.
The changeover task is activated at event point~$n+1$ in
Equation~(\ref{eqn:28a}) if it is a changeover from a task $i'\in
I^{\rm p1}$ to another task $i\in I^{\rm p1}$ or from a task
$i'\in I^{\rm np1}$ to any other task.

$$
wv(i'',n+1) = \sum_{i'\in I^{\rm p1} \cap I_u} \sum_{\ i\in I^{\rm
p1} \cap I_u, i'\neq i, Ctime_{i',i}>0} x(i',i,n) +
$$
\begin{equation}\label{eqn:28a}
\sum_{i'\in I^{\rm np1} \cap I_u} \sum_{\ i\in I_u, i'\neq i,
Ctime_{i',i}>0} x(i',i,n)  \ \ \forall u\in M, i''\in I^{\rm c}
\cap I_u, n<N^{\max}.
\end{equation}

The changeover task is activated at event point~$n+1$ in
Equation~(\ref{eqn:28b}) if it is a changeover from a task $i'\in
I^{\rm p1}$ to a task $i\in I^{\rm np1}$.

{\samepage
$$
wv(i'',n+1) = \sum_{i'\in I^{\rm p1} \cap I_u} \sum_{\  i\in
I^{\rm np1} \cap I_u, i'\neq i, Ctime_{i',i}>0} x(i',i,n)
$$
\begin{equation} \label{eqn:28b}
\forall u\in M, i''\in I^{\rm c1} \cap I_u, n<N^{\max}.
\end{equation}
}

The following two constraints correspond to equations (40c) and
(40d) in \cite{FloudasKTE2009} and express the duration of
changeover tasks:
$$
Tf(i'',n+1) - Ts(i'',n+1) = \sum_{i'\in I^{\rm p1}\cap I_u} \ \
\sum_{i\in I^{\rm p1}\cap I_u, i\neq i, Ctime_{i',i}>0}
Ctime_{i',i} \ x(i',i,n) +
$$
\begin{equation}
\sum_{i'\in I^{\rm np1}\cap I_u} \ \
 \sum_{i\in
I_u, i\neq i, Ctime_{i',i}>0} Ctime_{i',i} \ x(i',i,n)\ \forall
u\in M, i''\in I^{c} \cap I_u, n<N^{\max}.
\end{equation}

$$
Tf(i'',n+1) - Ts(i'',n+1) = \sum_{i'\in I^{\rm p1}\cap I_u} \ \
 \sum_{i\in I^{\rm np1}\cap I_u, i\ne i', Ctime_{i',i}>0}
 Ctime_{i',i} \ x(i',i,n)\
$$
\begin{equation}
\forall u\in M, i''\in I^{\rm c1} \cap I_u, n<N^{\max}.
\end{equation}

Note that (\ref{eqn:ch2}) and (\ref{eqn:ch3}) contain
$\Theta((N^{\max})^2 |I_u|^2)$ inequalities for each unit $u\in
M$, and each one has $\Theta(N^{\max} |I^{\rm p}\cap I_u| )$
summands. The number of binary variables $x(i',i,n)$ is
$\Theta(N^{\max} |I_u|^2)$. This is acceptable if the number of
tasks and event points is small.
For large number of tasks and event points, this creates serious
difficulties for the MIP solvers like CPLEX due to memory and CPU
time requirements. To overcome this problem, the more economical
model of changeover tasks was proposed as described in
Section~\ref{sec:new_mod}.

\subsection{Time Windows for Changeovers} \label{subsec:cb}

The changeover blockages between 4 a.m. and 7 a.m. and between
4~p.m. and 7~p.m. are modelled by explicit assignment of the event
points to appropriate time intervals. The number of event points
is proportional to the duration of the time interval. For example,
suppose that the horizon starts at 0:00 and finishes at 12:00 and
the total number of event points is $N^{\max} = 8$. The forbidden
interval for changeover tasks is from 4:00 to 7:00. Then there are
two allowed intervals from 0:00 to 4:00 and from 7:00 to 12:00
with the total duration of 9~hours. The number of event points
assigned to the first allowed interval  is  $[4\cdot 8/9] = [3.55]
= 4$, where $[\cdot]$ sign denotes rounding to the closest
integer. The other four event points are assigned to the second
interval.

\section{New Model of Changeover Tasks} \label{sec:new_mod}

Instead of binary variables $x(i',i,n)$ we use binary variables
$x(i,n)$ assuming that $x(i,n)$ equals 1 iff the task $i$ is the
first active task on its unit in the event points
$n,n+1,...,N^{\max}$. In other words, there is $n'\in
\{n,n+1,...,N^{\max}\}$ s.t. $wv(i,n')=1$ and all event points
$n,n+1,...,n'-1$ are empty on this unit.
Let us describe the constraints of the new model.\\

If a task $i$ is active at event point $n$ then $x(i,n)=1$:
\begin{equation}
\label{eq:new1} x(i,n) \ge wv(i,n)\ \forall u\in M, i\in I^{\rm
p}\cap I_u, n\le N^{\max}.
\end{equation}

If an event point $n$ contains no active tasks, then $x(i,n)$ is
``copied'' from event point $n+1$.
\begin{equation}
x(i,n) \ge x(i,n+1) - \sum_{i'\in I^{\rm p}\cap I_u,i'\neq i}
wv(i',n)\ \forall u\in M, i\in I^{\rm p}\cap I_u, n< N^{\max}.
\end{equation}

If task $i$ is followed by a different task $i'$ then a changeover
task must be activated at the event point $n+1$. At first we
consider the case $i, i' \in I^{\rm p1}$, so the changeover task
must be $i^{\rm c} \in I^{\rm c}$:
$$
wv(i^{\rm c},n+1) \ge \sum_{i'\in I^{\rm p1}\cap I_u,i'\neq i}
x(i',n+1) + wv(i,n) - 1
$$
\begin{equation}
\ \forall u\in M, i\in I^{\rm p1}\cap I_u, i^{\rm c}\in I^{\rm
c}\cap I_u, n< N^{\max}.
\end{equation}

In case $i\in I^{\rm np1}, i' \in I^{P}$,  the changeover task is
again $i^{\rm c} \in I^{\rm c}$:

$$
wv(i^{\rm c},n+1) \ge \sum_{i'\in I^{\rm p}\cap I_u,i'\neq i}
x(i',n+1) + wv(i,n) - 1
$$
\begin{equation}
\ \forall u\in M, i\in I^{\rm np1}\cap I_u, i^{\rm c}\in I^{\rm
c}\cap I_u, n< N^{\max}.
\end{equation}

In the last case $i\in I^{\rm p1}, i' \in I^{\rm np1}$,  the
changeover task  is $i^{\rm c1} \in I^{\rm c1}$:
$$
wv(i^{\rm c1},n+1) \ge \sum_{i'\in I^{\rm np1}\cap I_u} x(i',n+1)
+ wv(i,n) - 1
$$
\begin{equation}
\ \forall u\in M, i\in I^{\rm p1}\cap I_u, i^{\rm c1}\in I^{\rm
c1}\cap I_u, n< N^{\max}.
\end{equation}

If a changeover task is allocated then its duration is chosen as
$Ctime_{i,i'}$. We consider three cases again:
$$
Tf(i^{\rm c}, n+1) - Ts(i^{\rm c}, n+1) \ge \sum_{i'\in I^{\rm
p1}\cap I_u,i'\neq i}  Ctime_{i,i'} x(i',n+1)
 - H(1-wv(i,n))
$$
\begin{equation}
 \forall u\in M, i\in I^{\rm p1}\cap I_u, i^{\rm c}\in I^{\rm c}\cap I_u, n< N^{\max}.\ \end{equation}

$$
Tf(i^{\rm c}, n+1) - Ts(i^{\rm c}, n+1) \ge \sum_{i'\in I^{\rm
p}\cap I_u,i'\neq i} Ctime_{i,i'} x(i',n+1)
 - H(1-wv(i,n))
$$
\begin{equation}
 \forall u\in M, i\in I^{\rm np1}\cap I_u, i^{\rm c}\in I^{\rm c}\cap I_u, n< N^{\max}.\ \end{equation}

$$
Tf(i^{\rm c1}, n+1) - Ts(i^{\rm c1}, n+1) \ge \sum_{i'\in I^{\rm
np1}\cap I_u} Ctime_{i,i'} x(i',n+1)
 - H(1-wv(i,n))
$$
\begin{equation}
\label{eq:newLast}
 \forall u\in M, i\in I^{\rm p1}\cap I_u, i^{\rm c}\in I^{\rm c}\cap I_u, n < N^{\max}.\
 \end{equation}

The number of equations in (\ref{eq:new1})--(\ref{eq:newLast}) can
be estimated as~$\Theta(N^{\max} |I_u|)$ for each unit $u\in M$,
the number of binary variables $x(i,n)$ is~$\Theta(N^{\max}
|I_u|)$. This is by an order of magnitude smaller compared to
estimates in Section~\ref{sec:former_model} and allows for solving
the problem with larger horizons as in~\cite{BEFSK}.

Constraints (\ref{eq:new1})--(\ref{eq:newLast}) may be combined
with the changeover blockage mechanism~\cite{FloudasKTE2009}
described in Subsection~\ref{subsec:cb}. Nevertheless below we
propose a more flexible mechanism for blockage of changeover at
specific time windows which does not require any explicit
assignment of the event points to appropriate time intervals.

\subsection{Time Windows for Changeovers}

To model the changeover blockages between 4 a.m. and 7 a.m. and
between 4 p.m. and 7 p.m., we enumerate all the intervals where
changeovers are allowed and denote them by $[CT^s_k, CT^f_k], \
k\in K$. Here $K$ is the set of indices of changeover intervals.
For example, if the planning horizon starts at 0:00 a.m. then the
set of changeover intervals is $\{[0, 4], \ [7,16],\
[19,28],...\}$. Now we introduce new binary variables $z(i,k,n)$,
which equal~1 iff a changeover task $i$ at event point $n$ is
performed in the interval with number $k\in K$.

The following equations ensure that each changeover task is placed
in one time interval:
\begin{equation}
\sum_{k\in K} z(i,k,n) = wv(i,n) \ \forall i\in I^{\rm c}, n\leq
N^{\max}.
\end{equation}

The last pair of inequalities state that the changeover task is
performed in the chosen time interval:
\begin{equation}
Ts(i,n) \ge CT^s_k \ z(i,k,n) \ \forall i\in I^{\rm c},  k\in K,
n\leq N^{\max},
\end{equation}

\begin{equation}
Tf(i,n) \le CT^f_k + H(1-z(i,k,n))\ \forall i\in I^{\rm c},  k\in
K, n\leq N^{\max}.
\end{equation}

The blockages for the special changeovers from a task of $I^{\rm
p1}$ type to a task of $I^{\rm np1}$ type are modelled the same
way, but with a different set of time intervals.


\section{Conclusions}

A new formulation of changeover constraints for short-term
production scheduling problem is proposed. The new model requires
significantly less constraints compared to the original
formulation, which is important in case of large problem instances
where the memory requirements become a limiting factor for MIP
solvers.

\section{Acknowledgements} Partially supported by Russian Foundation
for Basic Research grants 12-01-00122 and 13-01-00862.

\end{document}